\magnification=1200
\baselineskip=16.7pt

\def\runningtitle{Approximation by conjugate integers}


\output={\ifnum\pageno>1
\headline{\hfill \rm \runningtitle \hfill }
\fi\plainoutput}

\def\and{\quad\hbox{and}\quad}

\def\Bbb{\bf}
\def\bA{{\bf A}}
\def\bC{{\Bbb C}}

\def\bQ{{\Bbb Q}}
\def\bR{{\Bbb R}}
\def\bZ{{\Bbb Z}}

\def\bfz{{\bf z}}

\def\cC{{\cal C}}
\def\cD{{\cal D}}
\def\cO{{\cal O}}
\def\cS{{\cal S}}

\def\deg{{\rm deg}}
\def\enum #1#2{#1_1,\dots,#1_{#2}}
\def\ep{\epsilon}

\def\Kbar{\bar{K}}
\def\Qbar{\bar{\bQ}}

\def\ndiv{\hbox{$\not|\,$}}
\def\proof{\medskip\noindent{\it Proof.}\quad }
\def\rank{{\rm rank}}
\def\Res{{\rm Res}}
\def\rmH{{\rm H}}
\def\section#1{\par\bigskip\goodbreak
  \hskip -1.2 true cm \noindent{\bf#1.\quad }}
\def\Vol{{\rm Vol}}

\def\ua{{\underline{a}}}
\def\uX{{\underline{X}}}
\def\uY{{\underline{Y}}}


\def\refBr{1}
\def\refBT{2}
\def\refBV{3}
\def\refARNAI{4}
\def\refCartan{5}
\def\refCassels{6}
\def\refCF{7}
\def\refFT{8}
\def\refLa{9}
\def\refLR{10}
\def\refMahler{11}
\def\refMahlerb{12}
\def\refMacfeat{13}
\def\refRo{14}
\def\refS{15}
\def\refTe{16}
\def\refWa{17}
\def\refWi{18}


%
%

\centerline{\bf Diophantine approximation by conjugate
algebraic integers}

\medskip
\centerline{\it Damien Roy\quad and\quad Michel Waldschmidt}

\quad
\hbox{\it \vtop{
 \hbox{D\'epartement de Math\'ematiques}
 \hbox{Universit\'e d'Ottawa, Canada}
}
\hfill\qquad \vtop{
 \hbox{Universit\'e P.~et M.~Curie (Paris VI)}
 \hbox{Institut de Math\'ematiques, Paris, France}
 } }

{\sevenrm\parindent=0pt
\footnote{}{2000 Mathematics classification. Primary 11J13;
Secondary 11J61, 11J82}
\footnote{}{Work of first author partially supported by
NSERC and CICMA}
}

\vskip1truecm
\noindent
{\bf Abstract:} Building on work of Davenport and Schmidt, we
mainly prove two results.  The first one is a version of
Gel'fond's transcendence criterion which provides a sufficient
condition for a complex or $p$-adic number $\xi$ to be algebraic
in terms of the existence of polynomials of bounded degree taking
small values at $\xi$ together with most of their derivatives. The
second one, which follows from this criterion by an argument of
duality, is a result of simultaneous approximation by conjugate
algebraic integers for a fixed number $\xi$ that is either
transcendental or algebraic of sufficiently large degree.  We also
present several constructions showing that these results are
essentially optimal.


\bigskip
\section{1. Introduction}
Motivated by work of Wirsing [\refWi], Davenport and Schmidt
investigated, in their 1969 seminal paper [\refARNAI], the
approximation of an arbitrary fixed real number $\xi$ by algebraic
integers of bounded degree.  They proved that if $n\ge 3$ is an
integer and if $\xi$ is not algebraic over $\bQ$ of degree at most
$(n-1)/2$ then there are infinitely many algebraic integers
$\alpha$ of degree at most $n$ such that
$$
|\xi-\alpha|\le c\rmH(\alpha)^{-[(n+1)/2]}
$$
where $c$ is a positive constant depending only on $n$ and $\xi$
and where $\rmH(\alpha)$ denotes the usual {\it height\/} of
$\alpha$, that is the maximum absolute value of the coefficients
of its irreducible polynomial over $\bZ$.  They also provided
refinements for $n\le 4$.  Recently, Bugeaud and Teuli\'e
revisited this result and showed in [\refBT] that we may also
impose that all approximations $\alpha$ have degree exactly $n$
over $\bQ$.  Moreover, a $p$-adic analog was proven by Teuli\'e
[\refTe].

Here we establish a similar result for simultaneous approximation
by several conjugate algebraic integers. In order to cover the
case where $\xi$ is a complex or $p$-adic number, we will assume
more generally that $\xi$ belongs to the completion of a number
field $K$ at some place $w$.

Thus, we fix an algebraic extension $K$ of $\bQ$ of finite degree
$d$. For each place $v$ of $K$, we denote by $K_v$ the completion
of $K$ at $v$ and by $d_v$ its local degree at $v$.  We also
normalize the corresponding absolute value $|\ |_v$ as in [\refBV]
by asking that, when $v$ is above a prime number $p$ of $\bQ$, we
have $|p|_v=p^{-d_v/d}$ and that, when $v$ is an Archimedean
place, we have $|x|_v=|x|^{d_v/d}$ for any $x\in\bQ$.  Then, our
result of approximation reads as follows:

\medskip
\proclaim Theorem A. Let $n$ and $t$ be integers with $1\le t\le
n/4$.  Let $w$ be a place of $K$ and let $\xi$ be an element of
$K_w$ which is not algebraic over $K$ of degree $\le (n+1)/(2t)$.
Assume further that $|\xi|_w\le 1$ if $w$ is ultrametric. Then
there exist infinitely many algebraic integers $\alpha$ which have
degree $n+1$ over $K$, degree $d(n+1)$ over $\bQ$ and admit, over
$K$, $t$ distinct conjugates $\enum \alpha t$ in $K_w$ satisfying
$$
\max_{1\le i\le t}|\xi-\alpha_i|_w
  \le c \rmH(\alpha)^{-(n+1)/(4dt^2)}
  \eqno(1.1)
$$
where $c$ is a constant depending only on $K$, $n$, $w$ and $\xi$.
\par

Note that, for $t=1$, $K=\bQ$ and $K_w=\bR$, this result is
comparable to Theorem 2 of [\refARNAI] mentioned above (with a
shift of $1$ in the degree of the approximation). Note also that,
if $w$ is ultrametric, any algebraic integer $\alpha$ in $K_w$
satisfies $|\alpha|_w\le 1$ and so the condition $|\xi|_w\le 1$ is
necessary to approximate $\xi$ by such numbers.

In section \S 10, we show that the exponent of approximation
$(n+1)/(4dt^2)$ in (1.1) is essentially best possible up to its
numerical factor of $1/4$ and that this factor cannot be replaced
by a real number greater than $2$, although its value can be
slightly improved using more precise estimates along the lines of
the present work.  For the sake of simplicity, we do not go into
such estimates here, nor do we try to sharpen the exponent of
approximation for small values of $n$.  It is difficult to predict
an optimal value for this exponent (see [\refRo]).

In answer to a question of K.~Tishchenko, we also show that one
cannot hope a similar exponent for simultaneous approximation of
$t$ numbers.  Taking $K=\bQ$ and $K_w=\bR$, we prove a result
which implies that, if $2\le t\le n$, then there exist a constant
$c>0$ and $t$ real numbers $\enum \xi t$ such that
$$
\max_{1\le  i\le t}|\xi_i-\alpha_i|
  \ge c \rmH(\alpha)^{-3n^{1/t}}
$$
for any choice of $t$ distinct conjugates $\enum \alpha t \in \bC$
of an algebraic number $\alpha$ of degree between $t$ and $n$.

The proof of Theorem A uses the same general strategy as Davenport
and Schmidt in [\refARNAI].  It relies on a duality argument
combined with the following version of Gel'fond's criterion of
algebraic independence where, for a polynomial $Q\in K[T]$, an
integer $j\ge 0$ and a
place $v$ of $K$, the notation $\|Q\|_v$
stands for the maximal
$v$-adic absolute value of the coefficients
of $Q$, while $Q^{(j)}$ denotes the $j$-th derivative of $Q$:

\proclaim Theorem B. Let $n$ and $t$ be integers with $1\le t \le
n/4$ and let $k=[n/4]$ denote the integral part of $n/4$.  Let $w$
be a place of $K$ and let $\xi$ be an element of $K_w$.  There
exists a constant $c>0$ which depends only on $K$, $n$, $w$ and $\xi$
and has the following property. Assume that, for each sufficiently
large real number $X$, there exists a nonzero polynomial $Q\in
K[T]$ of degree at most $n$ which satisfies $\|Q\|_v\le 1$ for
each place $v$ of $K$ distinct from $w$ and also
$$
\max_{0\le j\le n-t} |Q^{(j)}(\xi)|_w  \le cX^{-t/(k+1-t)}
  \and
\max_{n-t< j\le n} |Q^{(j)}(\xi)|_w  \le X. \eqno(1.2)
$$
Then, $\xi$ is algebraic over $K$ of degree $\le (n-k+1)/(2t)$.
\par

Note that Theorem B may still hold with an exponent smaller than
$t/([n/4]+1-t)$ in the above conditions (1.2).  However, we will
see in \S 3 that it would not hold with an exponent smaller than
$t/(n+1-t)$.

It is also interesting to compare this result with the criterion
of algebraic independence with multiplicities of [\refLR].  A main
difference is that the above theorem requires from the
polynomial $Q$ that a large proportion of its derivatives at
$\xi$ are small (at least three quarters of them), while the
conditions in Proposition 1 of [\refLR] are sharp only when a
small proportion of these derivatives are taken into account
(say, at most the first half of them).

This paper is organized as follows.  The next section \S2 sets the
various notions of heights that we use throughout the paper. The
results of duality which are needed to deduce Theorem A from
Theorem B are given in \S3, but the proof of this implication is
postponed to \S 9. Sections \S\S4--7 are devoted to preliminary
results towards the proof of Theorem B which is completed in \S 8.
In particular, we establish in \S4 a version of Gel'fond's
criterion (without multiplicities) which includes Theorem 2b of
[\refARNAI] and \S5 presents a height estimates which generalizes
Theorem 3 of [\refARNAI]. We conclude in \S 10 with two remarks on
the exponent of approximation in Theorem A.

\medskip
\noindent{\bf Notation.} Throughout this paper, $n$ denotes a
positive integer, $w$ denotes a place of $K$ and $\xi$ an element
of $K_w$.  For shortness, we will sometimes use the expressions
$a\ll b$ or $b \gg a$ to mean that given non-negative real numbers
$a$ and $b$ satisfy $a\le cb$ for some positive constant $c$ which
depends only on $K$, $n$, $w$ and $\xi$.  Overall, we tried to be
coherent with the notation of [\refARNAI].

%
%

\section{2. Heights}

Recall that $K$ is a fixed algebraic extension of $\bQ$ of degree
$d$. With the normalization of its absolute values given in the
introduction, the product formula reads
$$
\prod_v |a|_v = 1
$$
for any non-zero element $a$ of $K$.

Let $n$ be a positive integer.  For any place $v$ of $K$ and any
$n$-tuple $\ua=(\enum a n)\in K_v^n$, we define the {\it norm} of
$\ua$ as its maximum norm $\|\ua\|_v = \max\{|a_1|_v, \dots,
|a_n|_v\}$. Accordingly, the (absolute) {\it height} of a point
$\ua$ of $K^n$ is defined by
$$
\rmH(\ua) = \prod_v \|\ua\|_v.
$$

If $m$ is a positive integer with $1\le m\le n$ and if $M$ is an
$m\times n$ matrix with coefficients in $K_v$ for some place $v$
of $K$, we define $\|M\|_v$ as the norm of the ${n\choose
m}$-tuple formed by the minors of order $m$ of $M$ in some order.
When $M$ has coefficients in $K$, we define $\rmH(M)$ as the
height of the same point.  In particular, for an $m\times n$
matrix $M$ of rank $m$ with coefficients in $K$ we have
$\rmH(M)\ge 1$.

If $V$ is a subspace of $K^n$ of dimension $m\ge 1$, we define the
{\it height} $\rmH(V)$ of $V$ to be the height of a set of
Pl\"ucker coordinates of $V$.  In other words, we define
$\rmH(V)=\rmH(M)$ where $M$ is an $m\times n$ matrix whose rows
form a basis of $V$.  This is independent of the choice of $M$.
According to a well-known duality principle, we also have $\rmH(V)
= \rmH(N)$ where $N$ is any $(n-m)\times n$ matrix such that $V$
is the solution-set $\{ \ua\in K^n\;;\; N\ua=0 \}$ of the
homogeneous system attached to $N$ (see the formula (4), page 433
of [\refS]). When $V={0}$, we
set $\rmH(V)=1$.

We denote by $E_n$ the subspace of $K[T]$ consisting of all
polynomials of degree $\le n$, and for each place $v$ of $K$, we
denote by $E_{n,v}$ the closure of $E_n$ in $K_v[T]$.  We also
identify $E_n$ with $K^{n+1}$ and $E_{n,v}$ to $K_v^{n+1}$ by
mapping a polynomial $a_0+a_1T+\cdots+a_nT^n$ to the $(n+1)$-tuple
$(a_0,\dots,a_n)$ of its coefficients.  Accordingly, we define the
norm $\|P\|_v$ of a polynomial $P\in E_{n,v}$ as the maximum of
the absolute values of its coefficients and the height $\rmH(P)$
of a polynomial $P\in E_n$ as the height of the $(n+1)$-tuple of
its coefficients.  In the sequel, we will repeatedly use the fact
that, if $\enum P s\in K[T]$ have product $P=P_1\cdots P_s$ of
degree $\le n$, with $P\neq 0$ and $n\ge 1$, then we have
$$
e^{-n} \rmH(P_1)\cdots\rmH(P_s)
  < \rmH(P)
  < e^{n} \rmH(P_1)\cdots \rmH(P_s)
$$ as one gets for instance by comparing $ \|P_1\|_v\cdots
\|P_s\|_v$ and $\|P\|_v$ at all places $v$ of $K$ using the
various estimates of \S2 in Chapter 3 of [\refLa]. Finally, note
that, if $P$ is an irreducible polynomial of $K[T]$ of degree
$n$, if $\alpha$ is a root of $P$ in some extension of $K$ and if
$\deg(\alpha)$ denotes the degree of $\alpha$ over $\bQ$, then
there exist positive constants $c_1$ and $c_2$ depending only on
$n$ and $\deg(\alpha)$ such that
$$
c_1 H(\alpha)^n
   \le H(P)^{\deg(\alpha)}
   \le c_2 H(\alpha)^n.
$$
This follows from Proposition 2.5 in Chapter 3 of [\refLa]
applied once to $P$ and once to the irreducible polynomial of
$\alpha$ over $\bZ$ (since we defined $H(\alpha)$ to be the
height of the latter polynomial).

%
%

\section{3. Duality}

In this section, we fix a positive integer $n$, a place $w$ of $K$
and an element $\xi$ of $K_w$.  We define below a family of adelic
convex bodies and establish about them a result of duality that we
will need to deduce Theorem A from Theorem B.  We also describe
consequences of the adelic Minkowski's theorem of Macfeat
[\refMacfeat] and Bombieri and Vaaler [\refBV] for this type of
convex.

For any $(n+1)$-tuple of positive real numbers
$\uX=(X_0,X_1,\dots,X_n)$, we define an adelic convex body
$$
\cC(\uX)=\prod_v \cC_v(\uX)
  \subset \prod_v E_{n,v}
$$
by putting
$$
\cC_w(\uX)
  = \left\{ P\in E_{n,w} \,;\, \big|P^{(j)}(\xi)\big|_w \le X_j
           \hbox{ for $j=0,\dots,n$} \right\}
$$
at the selected place $w$ and
$$
\cC_v(\uX) = \{ P\in E_{n,v} \,;\, \|P\|_v \le 1 \}
$$
at the other places $v\neq w$.  For $i=1,\dots,n+1$, we denote by
$$
\lambda_i(\uX)=\lambda_i(\cC(\uX))
$$
the $i$-th minimum of $\cC(\uX)$ in $E_n$.  By definition, this is
the smallest positive real number $\lambda$ such that
$\lambda\cC(\uX)$ contains $i$ linearly independent elements of
$E_n$, where $\lambda\cC(\uX)$ is the adelic convex body whose
component at any Archimedean place $v$ consists of all products
$\lambda P$ with $P\in \cC_v(\uX)$ and whose component at any
ultrametric place $v$ is $\cC_v(\uX)$.

In order to apply the adelic Minkowski's theorem of [\refBV] in
this context, we identify each space $E_{n,v}$ with $K_v^{n+1}$ in
the natural way described in \S2. This identifies $\prod_v
E_{n,v}$ with $(K_\bA)^{n+1}$ where $K_\bA$ denotes the ring of
ad\`eles of $K$, and we use the same Haar measure as in [\refBV]
on this space. Explicitly, this means that, for an Archimedean
place $v$ of $K$, we choose the Haar measure on $K_v$ to be the
Lebesgue measure if $K_v=\bR$ and twice the Lebesgue measure if
$K_v=\bC$. For an ultrametric place $v$, we normalize the measure
so that the ring of integers $\cO_v$ of $K_v$ has volume
$|\cD_v|_v^{d/2}$ where $\cD_v$ denotes the local different of $K$
at $v$.  On each factor $E_{n,v}\simeq K_v^{n+1}$ we use the
product measure, and we take the product of these measures on
$\prod_v E_{n,v}$.

\medskip
\proclaim Lemma 3.1. There are two constants $c_1$ and $c_2$ which depend only
on $K$, $n$ and $w$ such that
$$
 c_1 (X_0\cdots X_n)^d\le \Vol(\cC(\uX)) \le c_2 (X_0\cdots X_n)^d
$$
for any $(n+1)$-tuple of positive real numbers $\uX =
(X_0,X_1,\dots,X_n)$.
\par

\proof Since the linear map from $E_{n,w}$ to itself sending a
polynomial $P(T)$ to $P(T+\xi)$ has determinant $1$, the volume of
$\cC_w(\uX)$ is equal to that of
$$
\left\{ P\in E_{n,w} \;;\; \big|P^{(j)}(0)\big|_w \le X_j
           \hbox{ for $j=0,\dots,n$} \right\}
\simeq
 \prod_{j=0}^n \big\{ a\in K_w \;;\; |j!a|_w \le X_j \big\}
$$
which in turn is bounded above and below by $c'_w (X_0\cdots X_n)^d$ and
$c''_w (X_0\cdots X_n)^d$  respectively, for some
positive constants $c'_w$ and $c''_w$  depending only on $K$, $n$ and $w$. For
the other places $v\neq w$, the volume of $\cC_v(\uX)$ is a
positive constant $c_v$ also depending only on $K$, $n$ and $v$,
with $c_v = 1$ for almost all places.  The conclusion follows.

\medskip
\proclaim Lemma 3.2. Let $\kappa$ be a positive real number and
let $\uX = (X_0,X_1,\dots,X_n)$ be a $(n+1)$-tuple of positive
real numbers with
$$
X_0\cdots X_n \ge c_1^{-1/d}(2/\kappa)^{n+1}
  \eqno{(3.1)}
$$
where $c_1$ is the constant of Lemma 3.1.  Then, $\kappa\cC(\uX)$
contains a non-zero element of $E_n$.
\par

\proof According to Theorem 3 of [\refBV], we have
$$
\big(\lambda_1(\uX)\cdots\lambda_{n+1}(\uX)\big)^d \Vol(\cC(\uX))
   \le 2^{d(n+1)}.
   \eqno{(3.2)}
$$
Since $\lambda_1(\uX)\le\cdots\le\lambda_{n+1}(\uX)$ and since
$\Vol(\cC(\uX)) \ge (2/\kappa)^{d(n+1)}$ by Lemma 3.1 and
condition (3.1), this implies $\lambda_1(\uX)\le \kappa$, as
required.

\medskip
Note that, for any integer $t$ with $1\le t\le n$ and any real
number $X\ge 1$, the condition (3.1) is satisfied with
$$
\kappa=1,
  \quad
X_0=\cdots=X_{n-t}=cX^{-t/(n+1-t)}
  \and
X_{n-t+1}=\dots=X_n=X
$$
for an appropriate constant $c$.  Then, the corresponding convex
body $\cC(\uX)$ contains a non-zero element of $E_n$.  In other
words, for any integer $t$ with $1\le t\le n$ and any real number
$X\ge 1$, there exists a nonzero polynomial $Q\in K[T]$ of degree
at most $n$ which satisfies $\|Q\|_v\le 1$ at each place $v$ of
$K$ distinct from $w$ and also
$$
\max_{0\le j\le n-t} |Q^{(j)}(\xi)|_w  \le cX^{-t/(n+1-t)}
  \and
\max_{n-t< j\le n} |Q^{(j)}(\xi)|_w  \le X.
$$
This justifies the remark made after the statement of Theorem B,
on comparing with the conditions (1.2) of this theorem.

\medskip
Our last objective of this section is to relate the successive
minima of a convex $\cC(X_0,\dots,X_n)$ with those of
$\cC(X_n^{-1},\dots,X_0^{-1})$.  We achieve this, following ideas
that go back to Mahler (see [\refMahler] and \S VIII.5 of
[\refCassels]), by showing that these convex bodies are almost
reciprocal with respect to some bilinear form $g$ on $E_n$. This
will require two lemmas.  The first one defines this bilinear form
$g$ and shows a translation invariance property of it.

\medskip
\proclaim Lemma 3.3.  Let $g\colon E_n\times E_n \to K$ be the
$K$-bilinear form given by the formula
$$
g(P,Q) = \sum_{j=0}^n (-1)^j P^{(j)}(0) Q^{(n-j)}(0)
$$
for any choice of polynomials $P,Q\in E_n$.  For each place $v$ of
$K$, denote by $g_v\colon E_{n,v}\times E_{n,v} \to K_v$ the
$K_v$-bilinear form which extends $g$.  Then, for any polynomials
$P,Q\in E_{n,w}$, we have
$$
g_w(P,Q) = \sum_{j=0}^n (-1)^j P^{(j)}(\xi) Q^{(n-j)}(\xi).
$$
\par

\proof For fixed $P, Q \in E_{n,w}$, the polynomial
$$
A(T) = \sum_{j=0}^n (-1)^j P^{(j)}(T) Q^{(n-j)}(T)
$$
has derivative
$$
A'(T)
  = \sum_{j=0}^{n-1} (-1)^j P^{(j+1)}(T) Q^{(n-j)}(T)
    + \sum_{j=1}^n (-1)^j P^{(j)}(T) Q^{(n-j+1)}(T)
  = 0.
$$
So $A(T)$ is a constant.  This implies $A(\xi)=A(0)=g_w(P,Q)$.

\medskip
Using this we get the following estimate.

\medskip
\proclaim Lemma 3.4. Let $\uX = (X_0,X_1,\dots,X_n)$ and $\uY =
(Y_0,Y_1,\dots,Y_n)$ be $(n+1)$-tuples of positive real numbers.
Suppose that, for each place $v$ of $K$, we are given polynomials
$P_v\in\cC_v(\uX)$ and $Q_v\in\cC_v(\uY)$.  Then, with the
notation of Lemma 3.3, we have
$$
\prod_v \big| g_v(P_v,Q_v) \big|_v
  \le (n+1)! \max_{0\le j\le n} X_j Y_{n-j}.
$$
\par

\proof For any place $v$ of $K$ with $v\neq w$, we have, if $v$ is
Archimedean,
$$
\eqalign{
\big| g_v(P_v,Q_v) \big|_v
 &\le (n+1)^{d_v/d}
      \max_{0\le j\le n} |P_v^{(j)}(0)|_v |Q_v^{(n-j)}(0)|_v \cr
 &\le \big((n+1)!\big)^{d_v/d} \|P_v\|_v \|Q_v\|_v \cr
 &\le \big((n+1)!\big)^{d_v/d}, \cr}
$$
and, if $v$ is ultrametric,
$$
\big| g_v(P_v,Q_v) \big|_v
  \le \max_{0\le j\le n} |P_v^{(j)}(0)|_v |Q_v^{(n-j)}(0)|_v
  \le \|P_v\|_v \|Q_v\|_v
  \le 1.
$$
Similarly, if $w$ is Archimedean, the formula of Lemma 3.3 leads
to
$$
\eqalign{
\big| g_w(P_w,Q_w) \big|_w
 &\le (n+1)^{d_w/d}
      \max_{0\le j\le n} |P_w^{(j)}(\xi)|_w |Q_w^{(n-j)}(\xi)|_w \cr
 &\le (n+1)^{d_w/d}
      \max_{0\le j\le n} X_j Y_{n-j} \cr}
$$
while, if $w$ is non-Archimedean, it gives
$$
\big| g_w(P_w,Q_w) \big|_w
  \le \max_{0\le j\le n} |P_w^{(j)}(\xi)|_w |Q_w^{(n-j)}(\xi)|_w
  \le \max_{0\le j\le n} X_j Y_{n-j}.
$$
The conclusion follows.

\medskip
\proclaim Proposition 3.5. Let $\uX=(X_0,\enum X n)$ be an
$(n+1)$-tuple of positive real numbers.  Define
$\uY=(Y_0,\dots,Y_n)$ where $Y_i=X_{n-i}^{-1}$ for $i=0,\dots,n$.
Then the products
$$
\lambda_i(\uX)\, \lambda_{n-i+2}(\uY), \quad (1\le i\le n+1),
$$
are bounded below and above by positive constants that depend only
on $K$, $n$ and $w$.
\par

\proof Fix an integer $i$ with $1\le i\le n+1$.  Put $\lambda =
\lambda_i(\uX)$ et $\mu = \lambda_{n-i+2}(\uY)$.  By definition of
the successive minima of a convex body, the polynomials of $K[T]$
contained in $\lambda \cC(\uX)$ generate a subspace $U$ of $E_n$
of dimension $\ge i$ while those contained in $\mu \cC(\uY)$
generate a subspace $V$ of $E_n$ of dimension $\ge n-i+2$. Since
the sum of these dimensions is strictly greater than that of $E_n$
and since the bilinear form $g$ of Lemma 3.3 is non-degenerate, it
follows that $U$ and $V$ are not orthogonal with respect to $g$.
Thus, there exist non-zero polynomials $P\in \lambda \cC(\uX)$ and
$Q\in \mu \cC(\uY)$ which belong to $E_n$ and satisfy $g(P,Q)\neq
0$. For any Archimedean place $v$ of $K$, we view $\lambda$ and
$\mu$ as elements of $K_v$ under the natural embedding of $\bR$ in
$K_v$ and define $P_v=\lambda^{-1}P$ and $Q_v=\mu^{-1}Q$.  For all
the other places of $K$, we put $P_v=P$ and $Q_v=Q$.  Then, we
have $P_v \in \cC_v(\uX)$ and $Q_v \in \cC_v(\uY)$ for all places
$v$ of $K$, and applying Lemma 3.4 we get
$$
\prod_{v\ndiv \infty} |g(P,Q)|_v \prod_{v| \infty}
\big|g_v(\lambda^{-1}P, \mu^{-1}Q)\big|_v
  \le (n+1)!.
$$
On noting that, for any Archimedean place $v$ of $K$, the real
numbers $\lambda$ and $\mu$ viewed as elements of $K_v$ satisfy
$|\lambda|_v=\lambda^{d_v/d}$ and $|\mu|_v=\mu^{d_v/d}$, we find
that the left hand side of the above inequality is
$$
\prod_{v|\infty} (\lambda\mu)^{-d_v/d} \prod_v |g(P,Q)|_v
  = (\lambda\mu)^{-1},
$$
by virtue of the product formula applied to the non-zero element
$g(P,Q)$ of $K$.  This shows $\lambda\mu\ge ((n+1)!)^{-1}$, and so
all products $\lambda_i(\uX) \lambda_{n-i+2}(\uY)$ are bounded
below by $((n+1)!)^{-1}$, for $i=1,\dots,n+1$.

On the other hand, applying Theorem 3 of [\refBV] to both
$\cC(\uX)$ and $\cC(\uY)$ (see (3.2) above), we find
$$
\eqalign{
\prod_{i=1}^{n+1} \big( \lambda_i(\uX) \lambda_{n-i+2}(\uY) \big)
  &= \left(\prod_{i=1}^{n+1} \lambda_i(\uX) \right)
     \left(\prod_{i=1}^{n+1} \lambda_i(\uY)\right) \cr
  &\le 4^{n+1} \Vol(\cC(\uX))^{-1/d}\, \Vol(\cC(\uY))^{-1/d} \cr
  & \le 4^{n+1} c_1^{-2/d} \cr}
$$
where $c_1$ is the constant of Lemma 3.1.  Thus the products
$\lambda_i(\uX) \lambda_{n-i+2}(\uY)$ are also bounded above by
$4^{n+1} c_1^{-2/d}((n+1)!)^n$.

%
%

\section{4. A version of Gel'fond's criterion}

Let $n$, $w$ and $\xi$ be as in the preceding section.  In this
section, we prove a specialized version of Gel'fond's
transcendence criterion which contains Theorem 2b of [\refARNAI]
and which we will need to conclude the proof of Theorem B.  It
applies as well to the situation of Lemma 12 in \S10 of
[\refARNAI].  For its proof, we need the following estimate (cf.\
Lemma 1 of [\refBr]).

\proclaim Lemma 4.1. Let $P,Q\in K[T]$ be relatively prime
non-zero polynomials of degree at most $n$. Then, we have
$$
1 \le c_3 \max\left\{ {|P(\xi)|_w\over \|P\|_w},
          {|Q(\xi)|_w\over \|Q\|_w} \right\}
          \rmH(P)^{\deg(Q)}
          \rmH(Q)^{\deg(P)}
$$
where $c_3=(2n)!$.
\par

\proof Since $P$ and $Q$ are relatively prime, their resultant
$\Res(P,Q)$ is a non-zero element of $K$.  For any place $v$ of
$K$, the usual representation of $\Res(P,Q)$ as a Sylvester
determinant leads to the estimate
$$
|\Res(P,Q)|_v \le c_v \|P\|_v^{\deg Q} \|Q\|_v^{\deg P}
$$
where $c_v=1$ if $v\ndiv\infty$ and $c_v=((2n)!)^{d_v/d}$ if
$v|\infty$.  Arguing as D.~Brownawell in the proof of Lemma 1 of
[\refBr], we also find
$$
|\Res(P,Q)|_w
  \le c_w \max\left\{ {|P(\xi)|_w\over \|P\|_w},
          {|Q(\xi)|_w\over \|Q\|_w} \right\}
          \|P\|_w^{\deg Q} \|Q\|_w^{\deg P}
$$
with the same value of $c_w$ as above.  The conclusion follows by
applying these estimates to the product formula $1 = \prod_v
|\Res(P,Q)|_v$.

\medskip
\proclaim Theorem 4.2.  Suppose that, for any sufficiently large
real number $X$, there is a non-zero polynomial $P=P_X\in K[T]$ of
degree $\le n$ and height $\le X$ such that
$$
{|P(\xi)|_w\over \|P\|_w} \le c_4^{-1} \rmH(P)^{-n} X^{-\deg(P)}
$$
where $c_4=e^{2n(n+1)}c_3^n$. Then, $\xi$ is algebraic over $K$ of
degree $\le n$ and the above polynomials vanish at $\xi$ for any
sufficiently large $X$.
\par

\proof We first reduce to a situation where we have monic
irreducible polynomials.  To this end, choose $X_0\ge 1$ such that
$P_X$ is defined for any $X\ge X_0$.  For a fixed $X\ge X_0$,
write $P=P_X$ in the form
$$
P = a Q_1\cdots Q_s
$$
where $a\in K^\times$ is the leading coefficient of $P$ and
$Q_1,\dots,Q_s$ are monic irreducible polynomials. Since $H(a)=1$
we have $H(Q_1)\cdots H(Q_s) \le e^n H(P)$ (see \S2) and so each
$Q_i$ has height at most $e^nX$.  Using this as well as the simple
estimate
$$
\|P\|_w
  \le |a|_w \prod_{i=1}^s \big( (1+\deg(Q_i)) \|Q_i\|_w \big)
  \le e^n |a|_w \prod_{i=1}^s \|Q_i\|_w,
$$
we deduce
$$
\prod_{i=1}^s
  \left(
  {|Q_i(\xi)|_w\over \|Q_i\|_w} \rmH(Q_i)^{n}
  \big( e^{n+1}X \big)^{\deg(Q_i)}
  \right)
\le e^{2n(n+1)} {|P(\xi)|_w\over \|P\|_w} \rmH(P)^{n} X^{\deg(P)}
\le c_3^{-n}.
$$
Writing $Y=e^nX$, we conclude that $P$ has at least one monic
irreducible factor $Q$ of degree $\le n$ and height $\le Y$ which
satisfies
$$
{|Q(\xi)|_w\over \|Q\|_w} \rmH(Q)^{n}
  \big( eY \big)^{\deg(Q)}
\le c_3^{-1}.
$$
Fix such a choice of polynomial $Q_Y=Q$ for each $Y\ge Y_0 =
e^nX_0$. Applying Lemma 4.1 to $Q_Y$ and $Q_{Y'}$ for values of
$Y$ and $Y'$ satisfying $e^nX_0\le Y\le Y' < eY$, we find that
these polynomials are not relatively prime.  Being monic and
irreducible, they are therefore equal to each other. So, we have
more generally that $Q_Y=Q_{Y_0}$ for any $Y\ge Y_0$.  As
$|Q_Y(\xi)|_w/\|Q_Y\|_w$ is bounded above by $c_3^{-1} (eY)^{-1}$
which tends to $0$ as $Y$ goes to infinity, this ratio must be $0$
independently of $Y\ge Y_0$.  This gives $Q_Y(\xi)=0$ for any
$Y\ge Y_0$ and therefore $P_X(\xi)=0$ for any $X\ge X_0$.

%
%

\section{5. A height estimate}

Here we establish a height estimate which, in our application,
will play the role of Theorem 3 of [\refARNAI].  Again, we start
with a lemma.

\proclaim Lemma 5.1. Let $\ell\ge 0$ be an integer and let
$x_0,\dots,x_\ell$ be indeterminates.  For any integer $k\ge 1$,
the set $\bZ[x_0,\dots,x_\ell]_k$ of homogeneous polynomials of
$\bZ[x_0,\dots,x_\ell]$ of degree $k$ is generated, as a
$\bZ$-module, by the minors of order $k$ of the $k\times (k+\ell)$
matrix
$$
\left.
R(k,\ell)
  = \pmatrix{
    x_0     &x_1    &\dots  &x_\ell    &0    &\dots  &0      \cr
    0       &x_0    &x_1    &\dots  &x_\ell  &\ddots &\vdots \cr
    \vdots  &\ddots &\ddots &\ddots &     &\ddots &0     \cr
    0       &\dots  &0      &x_0    &x_1  &\dots  &x_\ell   \cr}
\qquad
\right\} \quad \hbox{$k$ rows.}
$$

\proof We proceed by induction on $k+\ell$.  If $k=1$ or $\ell=0$
the result is clear.  Assume $k\ge 2$ and $\ell\ge 1$ and that the
result is true for a smaller number of rows or a smaller number of
indeterminates.  Denote by $M$ the subgroup of
$\bZ[x_0,\dots,x_\ell]_k$ generated by the minors of order $k$ of
the matrix $R(k,\ell)$.  The ring homomorphism $\varphi$ from
$\bZ[x_0,\dots,x_\ell]$ to $\bZ[x_0,\dots,x_{\ell-1}]$ sending
$x_\ell$ to $0$ and all other indeterminates to themselves maps
$M$ onto the subgroup of $\bZ[x_0,\dots,x_{\ell-1}]_k$ generated
by the minors of order $k$ of $R(k,\ell-1)$.  Thus, by the
induction hypothesis, we have
$$
\varphi(M) = \bZ[x_0,\dots,x_{\ell-1}]_k.
$$
On the other hand, the determinants of the $k\times k$
sub-matrices which contain the last column of $R(k,\ell)$ are the
products $x_\ell d$ where $d$ is a minor of order $k-1$ of
$R(k-1,\ell)$.  Thus, by the induction hypothesis, we also have
$$
M \supseteq x_\ell \bZ[x_0,\dots,x_\ell]_{k-1}
  = \bZ[x_0,\dots,x_\ell]_k \cap \ker(\varphi).
$$
These properties imply $M=\bZ[x_0,\dots,x_\ell]_k$.

\proclaim Proposition 5.2. Let $k$ and $\ell$ be integers with
$k\ge 1$ and $\ell\ge 0$.  For any $P\in E_\ell$, we have
$$
c^{-1} \rmH(P)^k \le \rmH(P\cdot E_{k-1}) \le c\rmH(P)^k
$$
where $c$ is a positive constant depending only on $k$ and $\ell$
and where $P\cdot E_{k-1}$ denotes the subspace of $E_{k+\ell-1}$
consisting of all products $PQ$ with $Q\in E_{k-1}$.
\par

\proof Write $P=a_0+a_1T+\cdots+a_\ell T^\ell$.  Then the height
of $P \cdot E_{k-1}$ is simply the height of the matrix with $k$
rows
$$
R = \pmatrix{
    a_0     &a_1    &\dots  &a_\ell &       &0 \cr
            &\ddots &\ddots &       &\ddots & \cr
    0       &       &a_0    &a_1    &\dots  &a_\ell \cr}.
$$
By virtue of the preceding lemma, every monomial of degree $k$
in $a_0,\dots,a_\ell$ can be expressed as a linear combination
of the minors of order $k$ of this matrix with integral
coefficients that do not depend on $P$.  Conversely, the
minors of order $k$ of $R$ can be written as linear
combinations of monomials of degree $k$ in $a_0,\dots,a_\ell$
with integral coefficients that do not depend on $P$.  Thus,
for each place $v$ of $K$, we have
$$
c_v^{-1} \|R\|_v \le \|P\|_v^k \le c_v \|R\|_v
$$
for some constant $c_v\ge 1$ independent of $P$, with $c_v=1$ when
$v$ is not Archimedean.  The conclusion follows with
$c=\prod_{v|\infty} c_v$.

%
%

\section{6. Construction of a polynomial}

Let $n$, $w$ and $\xi$ be as in \S 3.  We fix a non-decreasing
sequence of positive real numbers $X_0\le\cdots\le X_n$ and assume
that the corresponding convex body $\cC(X_0,\dots,X_n)$ contains a
non-zero polynomial $Q\in K[T]$.  Let
$$
V=\{P\in E_n\,;\, g(P,Q)=0 \}
$$
where $g\colon E_n\times E_n\to K$ is the $K$-bilinear form of
Lemma 3.3. For each integer $\ell$ with $0\le \ell\le n$, we
define a $K$-bilinear form $B_\ell\colon E_\ell\times E_{n-\ell}
\to K$ by the formula
$$
B_\ell(F,G)= g(FG,Q)
$$
for $F\in E_\ell$ and $G\in E_{n-\ell}$.  Its right kernel is
$$
V_\ell = \{ G\in E_{n-\ell} \,;\, G\cdot E_\ell \subseteq V \}
$$
We also denote by $B_{\ell,w}\colon E_{\ell,w}\times
E_{n-\ell,w}\to K_w$ the $K_w$-bilinear form which extends
$B_\ell$. Finally, we put
$$
y_i = (-1)^i i! Q^{(n-i)}(0)
  \and
z_i = (-1)^i i! Q^{(n-i)}(\xi),
  \quad
(0\le i\le n),
$$
and, for each integer $\ell=0,1,\dots,n$, we define
$$
M_\ell = \pmatrix{y_0 &y_1 &\dots &y_{n-\ell}\cr
                   y_1 &y_2 &\dots &y_{n-\ell+1}\cr
                \vdots &\vdots & &\vdots\cr
                   y_\ell &y_{\ell+1} &\dots &y_n\cr}
  \and
N_\ell = \pmatrix{z_0 &z_1 &\dots &z_{n-\ell}\cr
                   z_1 &z_2 &\dots &z_{n-\ell+1}\cr
                \vdots &\vdots & &\vdots\cr
                   z_\ell &z_{\ell+1} &\dots &z_n\cr}.
$$
With this notation, we will prove below a series of lemmas
leading, under some condition on $X_0,\dots,X_n$, to the
construction of a polynomial $P\in K[T]$ with several properties.
The method overall follows that of Davenport and Schmidt in
\S\S7--9 of [\refARNAI]. The first lemma is the following
observation:

\medskip
\proclaim Lemma 6.1.  Fix an integer $\ell$ with $0\le \ell\le n$.
Then,
\item{(i)} $M_\ell$ is the matrix of $B_\ell$ relative to the bases
  $\{1,T,\dots,T^\ell\}$ of $E_\ell$ and $\{1,T,\dots,T^{n-\ell}\}$ of
  $E_{n-\ell}$;
\item{(ii)} $N_\ell$ is the matrix of $B_{\ell,w}$ relative to
  the bases $\{1,T-\xi,\dots,(T-\xi)^\ell\}$ of $E_{\ell,w}$ and
  $\{1,T-\xi, \dots,(T-\xi)^{n-\ell}\}$ of $E_{n-\ell,w}$.
\par

\proof This follows upon noting that, for $i=0,\dots,\ell$ and
$j=0,\dots,n-\ell$, we have
$$
B_\ell(T^i,T^j) = g(T^{i+j},Q(T)) = y_{i+j}
$$
and, by Lemma 3.3,
$$
B_{\ell,w}((T-\xi)^i,(T-\xi)^j) = g_w((T-\xi)^{i+j},Q(T)) =
z_{i+j}.
$$

\medskip
In particular, this result implies that $M_\ell$ and $N_\ell$ have
the same rank for any value of $\ell$.  Note that the number of
rows of these matrices is less than or equal to their number of
columns if and only if $\ell\le n/2$.  Under this hypothesis, we
have the following estimates:

\medskip
\proclaim Lemma 6.2.  There are constants $c_5,c_6,c_7\ge 1$
depending only on $K$, $n$, $w$ and $\xi$ such that, for any
integer
$\ell$ with $0\le \ell \le n/2$, we have
\item{(i)} $\|N_\ell\|_w \le c_5 X_{n-\ell}\cdots X_n$,
\item{(ii)} $c_6^{-1}\|N_\ell\|_w \le \|M_\ell\|_w
  \le c_6\|N_\ell\|_w$,
\item{(iii)} $\rmH(M_\ell) \le c_7 \|N_\ell\|_w$ when $M_\ell$ has
  rank $\ell+1$.
\par

\proof (i) The upper bound on $\|N_\ell\|_w$ follows from the fact
that, for $i=1,\dots,\ell+1$, all the elements of the $i$-th row
of $N_\ell$ have their absolute value bounded above by a constant
times $X_{n+1-i}$.

(ii) By Lemma 6.1, $M_\ell$ and $N_\ell$ are matrices of
$B_{\ell,w}$ corresponding to different choices of bases for
$E_{\ell,w}$ and $E_{n-\ell,w}$.  Accordingly, we have
$$
M_\ell = {}^tU N_\ell V,
$$
where $U$ and $V$ are matrices of change of bases which depend
only on $\xi$, $\ell$ and $n$.  Since $U$ and $V$ are invertible,
this implies that any minor of order $\ell+1$ of $M_\ell$ (resp.\
$N_\ell$) can be expressed as a linear combination of the minors
of order $\ell+1$ of $N_\ell$ (resp.\ $M_\ell$) with coefficients
that are independent of $Q$, and the second assertion follows.

(iii) At any place $v$ of $K$ with $v\neq w$, the elements of
$M_\ell$ have their absolute value bounded above by a constant
which depends only on $n$ and which can be taken to be $1$ when
$v\ndiv\infty$.  So, the same is true of $\|M_\ell\|_v$.   The
height of $M_\ell$ is thus bounded above by a constant times
$\|M_\ell\|_w$ or, according to (ii), by a constant times
$\|N_\ell\|_w$.

\medskip
\proclaim Lemma 6.3. For any integer $\ell$ with $0\le \ell\le n$,
we have
$$
\dim V_\ell = n-\ell+1-\rank(M_\ell).
$$
When $M_\ell$ has rank $\ell+1$, we also have
$\rmH(V_\ell)=\rmH(M_\ell)$.
\par

\proof  A polynomial $P=a_0+a_1T+\cdots+a_{n-\ell}T^{n-\ell}$ of
$E_{n-\ell}$ belongs to $V_\ell$ if and only if for
$i=0,\dots,\ell$, it satisfies
$$
0 = g(T^iP(T),Q(T))
  = B_\ell(T^i,P(T))
  = \sum_{j=0}^{n-\ell} y_{i+j} a_j.
$$
Thus, identifying $E_{n-\ell}$ with $K^{n-\ell+1}$ in the usual
way, the subspace $V_\ell$ of $E_{n-\ell}$ is identified with the
solution space of the homogeneous system associated to $M_\ell$.
This proves the formula for $\dim V_\ell$.  Moreover, if $M_\ell$
has rank $\ell+1$, then, according to the duality principle
mentioned in \S2, we have $\rmH(V_\ell)=\rmH(M_\ell)$.

\medskip
\proclaim Lemma 6.4. Suppose that there exists an integer $h$ with
$1\le h\le n/2$ such that $M_{h-1}$ has rank $h$ and $M_h$ has
rank $\le h$.  Then, $V_{n-h}$ contains a non-zero element $P$.
Such a polynomial $P$ has degree $\le h$ and satisfies
$$
P\cdot E_{n-2h+1} = V_{h-1}.
$$
In particular, $P$ divides any polynomial of $V_{h-1}$.
\par

\proof  Since $M_{n-h}$ is the transpose of $M_h$, the two
matrices have the same rank.  By Lemma 6.3, this gives
$$
\dim V_{n-h} = (h+1)-\rank(M_{n-h}) \ge 1.
$$
So, $V_{n-h}$ contains a non-zero element $P$.  Using Lemma 6.3,
we also find
$$
\dim V_{h-1} =(n-h+2)-\rank(M_{h-1}) = n-2h+2.
$$
Since $V_{h-1}$ contains $P\cdot E_{n-2h+1}$, and since the
latter subspace of $E_{n-h+1}$ also has dimension $n-2h+2$,
this inclusion is an equality.

\medskip
\proclaim Lemma 6.5. Let $\ell$ and $t$ be integers with $0\le
\ell< n/2$ and $1\le t\le n-2\ell$.  Suppose that $N_\ell$ has
rank $\ell+1$ and that there exists a non-zero polynomial $P\in
K[T]$ such that $P\cdot E_{t-1} \subseteq V_\ell$ . Then, we have
$$
\left({|P(\xi)|_w\over \|P\|_w}\right)^t
  \ll {X_{n-t-\ell}\cdots X_{n-t} \over \|N_\ell\|_w}.
$$
\par

\proof  Denote by $\bfz_0,\dots,\bfz_{n-\ell}$ the columns of
$N_\ell$ and, for each integer $s$ with $1\le s\le t+1$, denote by
$N_\ell^{(s)}$ the sub-matrix of $N_\ell$ consisting of the
columns $\bfz_{s-1},\dots,\bfz_{n-\ell}$.  Observe that, since
$t\le n-2\ell$, these matrices all have at least $\ell+1$ columns.
Write
$$
P = b_0 + b_1 (T-\xi) + \cdots + b_h (T-\xi)^h
$$
where $h$ is the degree of $P$.  For any integer $s$ as above, we
have $(T-\xi)^{s-1}P(T)\in V_\ell$ and so, for $i=0,\dots,\ell$,
we find
$$
0 = B_{\ell,w}((T-\xi)^i,\,(T-\xi)^{s-1}P(T))
  = \sum_{j=0}^{h} z_{i+s-1+j} b_j.
$$
This means that the columns of $N_\ell$ satisfy the recurrence
relation
$$
b_0\bfz_{s-1}+b_1\bfz_{s}+\cdots+b_h\bfz_{s-1+h}=0,
  \quad
(1\le s\le t).
$$
Now, fix an integer $s$ with $1\le s\le t$ and choose indices
$j_0,j_1,\dots,j_\ell$ with $s-1\le j_0 < j_1 < \cdots < j_\ell
\le n-\ell$ such that
$$
\|N_\ell^{(s)}\|_w =
\big|\det(\bfz_{j_0},\dots,\bfz_{j_\ell})\big|_w. \eqno{(6.1)}
$$
If $j_0=s-1$, we find, using the recurrence relation,
$$
\eqalign{ |P(\xi)|_w \|N_\ell^{(s)}\|_w
  &= \big|\det(b_0\bfz_{s-1},\bfz_{j_1},\dots,\bfz_{j_\ell})
     \big|_w \cr
  &= \Big|\det\Big(-\sum_{j=1}^h b_j\bfz_{s-1+j},\bfz_{j_1},\dots,
         \bfz_{j_\ell}\Big)\Big|_w \cr
  &= \Big|\sum_{j=1}^h b_j \det(\bfz_{s-1+j},\bfz_{j_1},\dots,
         \bfz_{j_\ell})\Big|_w \cr
  &\le c \|P\|_w \|N_\ell^{(s+1)}\|_w \cr}
$$
for some positive constant $c$ depending only on $n$ and
$|\xi|_w$. If
$j_0\ge s$, this is still true because (6.1) then
implies
$\|N_\ell^{(s)}\|_w \le \|N_\ell^{(s+1)}\|_w$.  Since
$\|N_\ell^{(1)}\|_w = \|N_\ell\|_w \neq 0$, this inequality
implies by induction on $s$ that we have $\|N_\ell^{(s)}\|_w \neq
0$ for $s=1,\dots,t+1$, and therefore we can write
$$
{|P(\xi)|_w\over \|P\|_w}
  \le c {\|N_\ell^{(s+1)}\|_w \over \|N_\ell^{(s)}\|_w},
\quad (1\le s\le t).
$$
Multiplying term by term these inequalities, we get
$$
\left({|P(\xi)|_w\over \|P\|_w}\right)^t
  \le c^t {\|N_\ell^{(t+1)}\|_w \over \|N_\ell\|_w}
$$
and the conclusion follows upon noting that, for
$i=1,\dots,\ell+1$, the $i$-th row of $N_\ell^{(t+1)}$ has norm
$\ll X_{n-t-i+1}$ and thus $\|N_\ell^{(t+1)}\|_w \ll
X_{n-t-\ell}\cdots X_{n-t}$.

\medskip
\proclaim Proposition 6.6.  Let $k$ be an integer with $1\le k\le
n/2$.  Assume that there is an integer $t$ with $1\le t\le n+2-2k$
such that
$$
X_0\le \cdots\le X_{n-t} <1
  \and
1\le X_{n-t+1}\le \cdots\le X_n.
$$
Put $\delta = X_{n-t}$ and $Y = X_{n-t+1}\cdots X_n$, and assume
moreover that
$$
Y\delta^{k+1-t} < (c_5c_7)^{-1} \eqno(6.2)
$$
where $c_5$ and $c_7$ are defined in Lemma 6.2.  Then there exists
an integer $h$ with $1\le h\le k$ and a non-zero polynomial $P\in
K[T]$ of degree $\le h$ and height $\ll \delta^{-k/n}$ which
divides any polynomial of $V_{h-1}$ and satisfies
$$
\left({|P(\xi)|_w\over \|P\|_w}\right)^t
  \le c_8 \delta^h \rmH(P)^{-(n+2-2h)},
\eqno(6.3)
$$
where $c_8$ is a constant depending only on $K$, $n$, $w$ and $\xi$.
\par

\proof For any integer $\ell$ for which $M_\ell$ has rank
$\ell+1$, we find, using Lemma 6.2,
$$
\rmH(M_\ell)
  \le c_7 \|N_\ell\|_w
  \le c_5c_7 X_{n-\ell}\cdots X_n
  \le c_5c_7 Y\delta^{\ell+1-t}.
\eqno(6.4)
$$
Since we also have $\rmH(M_\ell)\ge 1$ for these values of $\ell$,
the assumption (6.2) implies that $M_k$ has rank $\le k$. The rank
of $M_0$ being 1, we conclude that there exists an integer $h$
with $1\le h\le k$ such that $M_{h-1}$ has rank $h$ and $M_h$ has
rank $\le h$. Then, according to Lemma 6.4, there exists a
non-zero polynomial $P\in E_h$ such that
$$
P\cdot E_{n-2h+1} = V_{h-1}.
$$
This implies that $P$ divides any polynomial of $V_{h-1}$ and, by
Proposition 5.2, that
$$
\rmH(P)^{n+2-2h} \ll \rmH(V_{h-1}) \ll \rmH(P)^{n+2-2h}.
\eqno(6.5)
$$
Combining Lemma 6.3 with (6.2) and (6.4) (for $\ell=h-1$), we also
find
$$
\rmH(V_{h-1})
  = \rmH(M_{h-1})
  \le c_5c_7 Y\delta^{h-t}
  < \delta^{-(k+1-h)}.
\eqno{(6.6)}
$$
Note that, since $k\le n/2$, the ratio $(k+1-h)/(n+2-2h)$ is a
decreasing function of $h$ in the range $1\le h\le k$. So, it is
bounded above by $k/n$.  Combining this observation with the above
estimates (6.5) and (6.6), we get
$$
\rmH(P)
  \ll \delta^{-(k+1-h)/(n+2-2h)}
  \ll \delta^{-k/n}. \eqno(6.7)
$$
Since $t\le n+2-2k$, we have $P\cdot E_{t-1} \subseteq P\cdot
E_{n-2h+1} \subseteq V_{h-1}$ and applying Lemma 6.5 with
$\ell=h-1$ gives
$$
\left({|P(\xi)|_w\over \|P\|_w}\right)^t
  \ll {\delta^h\over \|N_{h-1}\|_w}.
$$
Moreover, Lemma 6.2 (iii), Lemma 6.3 and (6.5) provide
$$
\|N_{h-1}\|_w
  \ge c_7^{-1} \rmH(M_{h-1})
   =  c_7^{-1} \rmH(V_{h-1})
   \gg \rmH(P)^{n+2-2h},
$$
and the conclusion follows.

\medskip
In our application, we will simply need the following consequence
of this proposition.

\medskip
\proclaim Corollary 6.7. Assume that all the hypotheses of
Proposition 6.6 are satisfied and that we have $\delta< c_8^{-1}$.
Then there exists an irreducible polynomial $P\in K[T]$ which
divides any polynomial of $V_{k-1}$ and satisfies
$$
\left({|P(\xi)|_w\over \|P\|_w}\right)^t
  \le c_9 \delta^{\deg P} \rmH(P)^{-(n+2-2k)}
\eqno(6.8)
$$
where $c_9=\max\{1,e^{n^2}c_8\}$.
\par

\proof Let $h$ and $P$ be as in the conclusion of the Proposition.
Since $\rmH(P)\ge 1$, the right hand side of (6.3) is bounded
above by $c_8\delta^h \le c_8\delta <1$.  So $P$ cannot be a
constant. Moreover, since $\deg(P)\le h\le k$, the same inequality
(6.3) gives
$$
\left({|P(\xi)|_w\over \|P\|_w}\right)^t
  \le c_8 \delta^{\deg(P)} \rmH(P)^{-(n+2-2k)}.
$$
Write $P$ as a product
$$
P= P_1\cdots P_s
$$
of irreducible polynomials of $K[T]$.  Then the above inequality
leads to
$$
\prod_{i=1}^s
  \left(
  \left({|P_i(\xi)|_w\over \|P_i\|_w}\right)^t
        \delta^{-\deg(P_i)} \rmH(P_i)^{n+2-2k}
  \right)
  \le e^{n^2}c_8 \le c_9^s.
$$
So, a least one factor of the product on the left must be bounded
above by $c_9$.  The corresponding polynomial $P_i$ divides every
element of $V_{k-1}$ since it divides $P$ and $V_{k-1}$ is
contained in $V_{h-1}$.

\medskip
Note that this statement provides no upper bound on the degree and
height of $P$.  We will get such upper bounds by an indirect
argument, using the construction of an auxiliary polynomial in the
next section.

%
%

\section{7. Degree and height estimates}

The notation is as in the preceding section.  We assume that the
adelic convex body $\cC(X_0,\dots,X_n)$ contains a non-zero
polynomial $Q$ of $K[T]$ for some non-decreasing sequence of
positive real numbers $X_0\le\cdots\le X_n$, and we define
corresponding subspaces $V_\ell$ of $E_{n-\ell}$ for
$\ell=0,\dots,n$ as in \S 6.

\proclaim Lemma 7.1. Put $c=((n+1)!)^{-2}$. Then, for any integer
$\ell$ with $0\le \ell \le n$ we have
$$
\cC(cX_n^{-1},\dots,cX_{\ell}^{-1}) \cap E_{n-\ell}
  \subseteq V_{\ell}.
\eqno(7.1)
$$
\par

\proof Let $\ell$ be an integer with $0\le \ell\le n$, and let $G$
be an element of the left hand side of (7.1).  We need to show
that $g(T^m G,Q)=0$ for $m=0,\dots,\ell$.  To this end, we
proceed by induction.  We fix an integer $m$ with $0\le m\le
\ell$ and assume, when $m\ge 1$, that we have $g(T^jG,Q)=0$ for
$j=0,\dots,m-1$.  Let $P=T^mG(T)$.  We define $P_w=(T-\xi)^m
G(T)$ and $Q_w=Q$ and, for the other places $v\neq w$ of $K$, we
put $P_v=P$ and $Q_v=Q$.  These polynomials satisfy
$P_v\in\cC_v(\uY)$ and $Q_v\in\cC_v(\uX)$ for each place $v$ of
$K$, where $\uY=(Y_0,\dots,Y_n)$ denotes the $(n+1)$-tuple of
positive real numbers given by $Y_i = n! c X_{n-i}^{-1}$ for
$i=0,\dots,n$. Moreover, the hypotheses imply
$$
g_w(P_w,Q_w) = g(P,Q)
$$
since the difference $P_w-P$ can be written as a linear
combination of $G,\dots,T^{m-1}G$ with coefficients in $K_w$ in
the case $m\ge 1$ and is zero when $m=0$.  Using Lemma 3.4, we
therefore get
$$
\prod_v |g(P,Q)|_v
  = \prod_v |g_v(P_v,Q_v)|_v
  \le (n+1)! n! c
  < 1.
$$
By the product formula, this implies $g(P,Q)=0$.

\medskip
\proclaim Proposition 7.2.  There is a constant $c_{10}>0$ which
depends only on $K$, $n$, $w$ and $\xi$ and has the following property.
Suppose that $\ell$ and $u$ are non-negative integers with $\ell +
u < n$, such that
$$
X_n^{u+1}X_{n-1}\cdots X_{\ell+u} \le c_{10}. \eqno(7.2)
$$
Then, there is a non-zero polynomial $G$ of $K[T]$ of degree $\le
n-\ell$ and height $\ll X_{\ell+u}^{-1}$ such that $G^{(i)}\in
V_{\ell}$ for $i=0,\dots,u$.
\par

\proof Let $c$ be as in Lemma 7.1. Put $\kappa=(n!)^{-1}$ and
define real numbers $Y_0,\dots,Y_{n-\ell}$ by
$$
Y_i = \cases{ cX_n^{-1} &for $i=0,\dots,u$,\cr\cr
              cX_{n-i+u}^{-1} &for $i=u+1,\dots,n-\ell$.\cr}
$$
Lemma 3.2 shows that the convex $\kappa\cC(Y_0,\dots,Y_{n-\ell})$
contains a non-zero element $G$ of $E_{n-\ell}$ if the condition
(7.2) is satisfied for a sufficiently small constant $c_{10}>0$. Such
a polynomial has height $\rmH(G) \ll Y_{n-\ell}\ll
X_{\ell+u}^{-1}$. Moreover, for $i=0,\dots,u$, we find
$$
G^{(i)}
  \in \cC(Y_i,\dots,Y_{n-\ell})
  \subseteq \cC(cX_n^{-1},\dots,cX_{\ell}^{-1})
$$
and so $G^{(i)}\in V_{\ell}$ by Lemma 7.1.

\medskip
We will apply this proposition in the following context:

\medskip
\proclaim Corollary 7.3. Let $\ell$ and $u$ be as in Proposition
7.2, and assume that there exists an irreducible polynomial $P\in
K[T]$ which divides every element of $V_\ell$.  Then, we have
$$
\deg(P) \le {n-\ell\over u+1}
  \and
\rmH(P) \ll X_{\ell+u}^{-1/(u+1)}.
$$
\par

\proof  The hypotheses imply that $P$ divides all derivatives of
the polynomial $G$ of Proposition 7.2, up to order $u$.  So,
$P^{u+1}$ divides $G$ and the conclusion follows.

%
%

\section{8. Proof of Theorem B}

Let the notation be as in Theorem B and assume that the hypothesis
of this theorem holds with a constant $c<\min\{1,(c_5c_7)^{-1}\}$.
Then, for any real number $X\ge 1$ the condition (6.2) of
Proposition 6.6 is satisfied with
$$
\delta=X_0=X_1=\cdots=X_{n-t}=cX^{-t/(k+1-t)}
  \and
X_{n-t+1}=\cdots=X_n=X.
$$
Moreover, if $X$ is sufficiently large, the hypothesis of Theorem
B is that the corresponding convex $\cC(\uX)$ with
$\uX=(X_0,\dots,X_n)$ contains a non-zero element $Q$ of $E_n$.
Since $t\le n+2-2k$, we may then apply Corollary 6.7. It shows
that, if $X$ is sufficiently large so that $\delta < c_8^{-1}$,
then there is an irreducible polynomial $P\in K[T]$ which divides
every element of the vector space $V_{k-1}$ attached to $Q$ and
satisfies
$$
{|P(\xi)|_w\over \|P\|_w}
  \ll \rmH(P)^{-(n+2-2k)/t} \delta^{\deg(P)/t}.
$$
Since $c\le 1$ and $n-t \ge 2t+k-2$, we also find
$$
X_n^{2t}X_{n-1}\cdots X_{2t+k-2}
   =  X^{3t-1} \delta^{n-3t-k+3}
  \le X^{3t-1} \delta^{3k-3t+3}
  \le X^{-1}.
$$ So, the condition (7.2) of Proposition 7.2 is satisfied with
$\ell=k-1$ and $u=2t-1$ provided that $X$ is sufficiently large.
Assuming that this is the case, Corollary 7.3 then shows $$
\deg(P) \le {n-k+1 \over 2t}
  \and
\rmH(P) \le \kappa \delta^{-1/(2t)}
$$
for some constant $\kappa>0$.  Putting $m=[(n-k+1)/(2t)]$ and
$Y=\kappa \delta^{-1/(2t)}$, and noting that $(n+2-2k)/t \ge m$,
we thus have found the existence of a polynomial $P\in K[T]$ of
degree $\le m$ and height $\le Y$ such that
$$
{|P(\xi)|_w\over \|P\|_w}
  \ll \rmH(P)^{-m} Y^{-2\deg(P)}.
$$
Since $Y$ is a monotone increasing unbounded continuous function
of $X$, for $X\ge 1$, Theorem 4.2 then shows that $\xi$ is
algebraic over $K$ of degree $\le m$.

%
%
\section{9. Proof of Theorem A}

We first generalize the construction of Davenport and Schmidt in
\S2 of [\refARNAI].

\medskip \proclaim Lemma 9.1. Let $t$ be an integer with $1\le
t\le n$ and let $\delta$, $\kappa$ and $Y$ be real numbers with
$0 < \delta < 1 < Y$ and $\kappa \ge 1$. Assume that
$\uY=(Y_0,\dots,Y_n)$ is a $(n+1)$-tuple of positive real numbers
satisfying $$ \lambda_{n+1}(\uY) \le \kappa
  \and
Y_j\le \cases{
 Y\delta^{t-j} &for $j=0,\dots,t-1$,\cr
 Y &for $j=t,\dots,n.$\cr}
$$ Assume moreover that $|\xi|_w\le 1$ in the case where $w$ is
ultrametric. Then there exists a monic polynomial $P\in \cO_K[T]$
which is irreducible over $K$ of degree $n+1$, has height $H(P)\ll
\kappa Y$ (as defined in \S2), admits $d$ distinct conjugates over
$\bQ$, and has at least $t$ distinct roots in the closed disk of
$K_w$ centered at $\xi$ of radius $\delta$. \par

\proof Let $\epsilon$ be a fixed but arbitrarily small positive
real number with $\epsilon \le 1$. Put $\delta_0 =
\min\{\delta,\epsilon\}$ and choose elements $P_1,\dots,P_{n+1}$
in $E_n$ realizing the successive minima of $\cC(\uY)$ in $E_n$.
Since $\lambda_{n+1}(\uY)\le \kappa$, these polynomials all
belong to $\kappa\cC(\uY)$.  In particular, they have integral
coefficients at any ultrametric place $v$ of $K$ distinct from
$w$.  Moreover, they form a basis of $E_n$ over $K$.  We will
construct the required polynomial $P$ in the form
$$
P(T) = T^{n+1} + \sum_{i=1}^{n+1} b_i P_i(T)
$$
for suitable elements $\enum b {n+1}$ of $K$. Each $b_i$ will
be obtained as solution of a system of inhomogeneous inequalities
using the strong approximation theorem (see Theorem 3, page 440
of [\refMahlerb] or \S15 of [\refCF]). According to this result,
there is a constant $C>0$ depending only on $K$ with the
following property. For any finite set $\cS$ of places of $K$,
any choice of elements $\beta_v\in K_v$, ($v\in\cS$), and any
choice of positive real numbers $\epsilon_v$, ($v\in \cS$), with
$\prod_{v\in\cS} \epsilon_v \ge C$, there exists an element $b$
of $K$, satisfying $|b-\beta_v|_v \le \epsilon_v$ for $v\in \cS$,
and $|b|_v\le 1$ for $v\notin \cS$.

To ensure that $P$ is irreducible over $K$ and admits $d$ distinct
conjugates over $\bQ$, we proceed essentially as Bugeaud and
Teuli\'e in [\refBT].  We choose a prime number $q$ of $\bZ$ which
splits completely in $\cO_K$ into a product of $d$ distinct prime
ideals none of which defines the place $w$.  We fix a place $v_0$
among the corresponding $d$ places of $K$ above $q$ and we choose
an element $\pi$ of $K$ satisfying $|\pi|_{v_0}=|q|_{v_0}$ and
$|\pi|_v=1$ for $v|q$ with $v\neq v_0$. We write
$$
\pi =
  \sum_{i=1}^{n+1} \gamma_i P_i(T)
$$
with $\enum \gamma {n+1} \in K$ and we ask that
$$
|b_i-\gamma_i|_v \le |q|_v^2 = q^{-2/d},
  \quad (v|q,\ 1\le i\le n+1).
\eqno(9.1)
$$
Under these conditions, the corresponding polynomial $P$ satisfies
$$
\| P(T) - T^{n+1} - \pi \|_v
  = \Big\| \sum_{i=1}^{n+1} (b_i-\gamma_i) P_i(T) \Big\|_v
  \le |q|_v^2
$$ for $v|q$.  Thus, $P$ has integral coefficients at the places
of $K$ above $q$.  Since $\pi$ is a uniformising  parameter for
$v_0$, the above relation implies that $P$ is an Eisenstein
polynomial of $K[T]$ at $v_0$ and thus it is irreducible over $K$
(see for instance Theorem 24 in \S~3, Chapter III of [\refFT]).
Moreover this relation also gives $|P(0)|_{v_0}<1$ and
$|P(0)|_v=1$ for $v|q$ with $v\neq v_0$.  Thus the constant
coefficient $P(0)$ of $P$ admits $d$ distinct conjugates over
$\bQ$.

To ensure that $P$ has $t$ roots close to $\xi$ in $K_w$, we fix a
monic polynomial $B\in K_w[T]$ of degree $t$ with $t$ (simple)
distinct roots in the open unit disk $D_w=\{z\in K_w\; ; \;
|z|_w<1\}$ of $K_w$ and we use the fact that, by explicit forms of
the inverse function theorem such as Theorem 4.4.1 in Chapter I of
[\refCartan], any polynomial $S(T)\in K_w[T]$ of degree $\le n+1$
for which $\|S-B\|_w$ is sufficiently small also has $t$ distinct
roots in $D_w$.  We proceed as follows:

If $w$ is ultrametric, lying above an ordinary prime number $p$,
we choose an element $r$ of $K_w$ with $p^{-1}\delta_0\le |r|_w
\le \delta_0$ and put $s=r^t$. If $w$ is Archimedean, we choose
$r,s\in K_w$ with $|r|_w=\delta_0$ and $|s|_w = \kappa^{d_w/d}
\epsilon^{-t-2} \delta_0^t Y$. In both cases, we define
$$
R(T) =
  (T-\xi)^{n+1} + s B\left({T-\xi\over r}\right) \in K_w[T].
$$
We write this polynomial in the form
$$
R(T) =
  T^{n+1} +  \sum_{i=1}^{n+1} \theta_i P_i(T)
$$
with $\enum \theta {n+1} \in
K_w$ and ask that
$$
|b_i-\theta_i|_w
  \le
\cases{
     \epsilon^{-1}   &if $w$ is Archimedean,\cr
     \epsilon^{n+1} Y^{-1} &if $w$ is ultrametric. \cr}
  \eqno{(9.2)}
$$
The polynomial $S=s^{-1}P(rT+\xi)$ then satisfies
$$
\| S - B \|_w
  = \left\|
    s^{-1}r^{n+1}T^{n+1}
    + s^{-1} \sum_{i=1}^{n+1} (b_i-\theta_i) P_i(rT+\xi)
    \right\|_w
  \ll \epsilon
$$
using $|s|_w^{-1}|r|_w^{n+1} \ll \delta_0 \ll \epsilon$ and noting
that, for $i=1,\dots,n+1$, we have
$$
\|P_i(rT+\xi)\|_w
  \ll \cases{
  \kappa^{d_w/d}\delta^t Y &if $w$ is Archimedean,\cr
  \delta^t Y               &if $w$ is ultrametric.\cr}
$$
If $\epsilon$ is sufficiently small, this implies that $S$ has $t$
roots in the disk $D_w$ and therefore, that $P$ has at least $t$
distinct roots in the disk of $K_w$ centered at $\xi$ with radius
$\delta$.

If $w$ is Archimedean and again if $\epsilon$ is small enough, the
strong approximation Theorem allows us to require, aside from
(9.1) and (9.2), that
$$
|b_i|_v \le 1, \quad (1\le i\le n+1),
  \eqno{(9.3)}
$$
for all places $v$ of $K$ with $v\neq w$ and $v\ndiv q$.  Then $P$
has integral coefficients at $v$ for each ultrametric place $v$ of
$K$ and therefore it has coefficients in $\cO_K$.  Moreover, as
we may take $\ep \gg 1$, we find $\|P\|_w \ll \kappa^{d_w/d}
\ep^{-t-2} Y \ll \kappa^{d_w/d} Y$.  Since $\|P\|_v \ll
\kappa^{d_v/d}$ for all the other Archimedean places $v$ of $K$,
this implies $\rmH(P)\ll \prod_{v|\infty} \|P\|_v \ll \kappa Y$.

If $w$ is ultrametric, we choose an Archimedean place $v_1$.  We
require that (9.1) and (9.2) hold, that (9.3) holds for all places
$v$ of $K$ with $v\neq w$, $v\neq v_1$ and $v\ndiv q$, and that
$$
|b_i|_{v_1} \le \epsilon^{-n-2}Y, \quad (1\le i\le n+1).
$$
Again, the strong approximation Theorem shows that these
conditions have solutions $b_i\in K$ for $i=1,\dots,n+1$ provided
that $\epsilon$ is small enough.  Then, the corresponding
polynomial $P$ has integral coefficients at $v$ for each
ultrametric place $v$ of $K$ with $v\neq w$.  At the place $w$,
we find
$$
\| P - R \|_w
  = \left\| \sum_{i=1}^{n+1} (b_i-\theta_i) P_i(T) \right\|_w
  \le \epsilon^{n+1} Y^{-1} \max_{1\le i\le n+1} \|P_i\|_w
  \ll \epsilon^{n+1}.
$$
Moreover, $R$ has coefficients in $\cO_w$ since $|\xi|_w\le 1$
and $\|B\|_w\le 1$.  Thus, $P$ also has coefficients in $\cO_w$
if $\epsilon$ is sufficiently small, and then it has
coefficients in $\cO_K$.  Since we may take $\ep \gg 1$, this
gives $\|P\|_{v_1} \ll
\kappa^{d_{v_1}/d} Y$ and, since $\|P\|_v
\ll \kappa^{d_v/d}$ for
all Archimedean places $v\neq v_1$, we
find $\rmH(P)\ll
\prod_{v|\infty} \|P\|_v \ll \kappa Y$.

\bigskip
\noindent
{\bf Proof of Theorem A.}
\par
Let $t$, $n$ and $\xi$ be as in Theorem A, let $k=[n/4]$, and let
$c$ be the constant of Theorem B corresponding to these data.
Since $\xi$ is not algebraic over $K$ of degree $\le
(n-k+1)/(2t)$, Theorem B shows that there are arbitrary large
positive real numbers $X$ for which the $(n+1)$-tuple
$\uX=(X_0,\dots,X_n)$ given by
$$
X_j= \cases{
 cX^{-t/(k+1-t)} &for $j=0,\dots,n-t$,\cr
 X &for $j=n-t+1,\dots,n$\cr}
$$
satisfies $\lambda_1(\uX)>1$.  According to Proposition 3.5, this
implies that the $(n+1)$-tuple $\uY=(Y_0,\dots,Y_n)$ given by
$$
Y_j = X_{n-j}^{-1} = \cases{
 X^{-1} &for $j=0,\dots,t-1$,\cr
 c^{-1}X^{t/(k+1-t)} &for $j=t,\dots,n$\cr}
$$
satisfies $\lambda_{n+1}(\uY) \le \kappa$ with a constant
$\kappa\ge 1$ depending only on $K$, $n$ and $w$.  Assuming $X$
sufficiently large, we may thus apply Lemma 9.1 with
$$
Y=c^{-1}X^{t/(k+1-t)} \and \delta=c^{1/t}X^{-(k+1)/(t(k+1-t))}.
$$
It shows the existence of a monic polynomial $P\in \cO_K[T]$ which
is irreducible over $K$ of degree $n+1$ and height $\ll \kappa Y$,
admits $d$ distinct conjugates over $\bQ$, and has at least $t$
distinct roots $\enum \alpha t \in K_w$ with
$$
\max_{1\le i\le t} |\xi-\alpha_i|_w
  \ll Y^{-(k+1)/t^2}.
\eqno{(9.4)}
$$
In particular, $\alpha=\alpha_1$ is an algebraic integer of degree
$n+1$ over $K$ and degree $d(n+1)$ over $\bQ$.  From the remark at
the end of \S2, we find $\rmH(\alpha) \ll \rmH(P)^d \ll Y^d$.
Combining this with (9.4), we obtain that the conjugates $\enum
\alpha t$ of $\alpha$ over $K$ satisfy
$$
\max_{1\le i\le t} |\xi-\alpha_i|_w
  \ll \rmH(\alpha)^{-(k+1)/(dt^2)}
  \ll \rmH(\alpha)^{-(n+1)/(4dt^2)}.
$$
Moreover, since the right hand side of (9.4) can be made
arbitrarily small by choosing $X$ sufficiently large, we produce
an infinity of such numbers $\alpha$ by varying $X$.

%
%
\section{10. Remarks on simultaneous approximations}

We fix a place $w$ of $K$ and an algebraic closure $\Kbar_w$ of
$K_w$, and we extend the absolute value $|\ |_w$ of $K_w$ to an
absolute value of $\Kbar_w$ also denoted $|\ |_w$.  Our first
result below shows that, for $t\ge 2$, the exponent
$(n+1)/(4dt^2)$ in the inequality (1.1) of Theorem A cannot be
replaced by a real number greater than $2n/(dt(t-1))$.

\proclaim Proposition 10.1. Let $n$ and $t$ be integers with $2\le
t\le n$, and let $\xi$ be an element of $K_w$.  There exists a
constant $c=c(n,t)>0$ such that, for any algebraic number $\alpha$
of degree $n$ over $K$ and any choice of $t$ distinct conjugates
$\enum \alpha t$ of $\alpha$ in $\Kbar_w$, we have
$$
\max_{1\le i\le t} |\xi-\alpha_i|_w
  \ge c\rmH(\alpha)^{-\textstyle{2n(n-1)\over \deg(\alpha)t(t-1)}}
$$
where $\deg(\alpha)$ denotes the degree of $\alpha$ over $\bQ$.
\par

\proof  Let $P\in K[T]$ be an irreducible polynomial of degree
$n$, let $a_0$ be its leading coefficient, and let $\enum \alpha
n$ be the roots of $P$ ordered so that
$$
|\xi-\alpha_1|_w\le\cdots\le|\xi-\alpha_n|_w.
$$
The discriminant of $P$ is the non-zero element of $K$ given by
$$
{\rm Disc}(P)
  = {a_0}^{2(n-1)} \prod_{1\le i<j \le n}(\alpha_i-\alpha_j)^2.
$$
Using the estimates
$$
|\alpha_i-\alpha_j|_w
  \le \cases{
  2^{c_w} \max\{1,|\alpha_i|_w\}\max\{ 1,|\alpha_j|_w\}|\xi-\alpha_j|_w
  &when $1\le i<j\le t$, \cr\cr
  2^{c_w} \max\{1,|\alpha_i|_w\}\max\{ 1,|\alpha_j|_w\}
  &otherwise, \cr}
$$
with $c_w=0$ if $w$ is ultrametric and $c_w=d_w/d$ otherwise we
find
$$
|{\rm Disc}(P)|_w
  \le 2^{n(n-1)c_w} M_w(P)^{2(n-1)}
\prod_{j=1}^t|\xi-\alpha_j|_w^{2(j-1)},
$$
where $M_w(P) = |a_0|_w\prod_{i=1}^m \max\{1,|\alpha_i|_w\}$
denotes the Mahler measure of $P$ at $w$.  Since $M_w(P) \le
(n+1)^{c_w/2} \|P\|_w$ (see Chapter 3 of [\refLa]), this gives
$$
|{\rm Disc}(P)|_w
  \le \big(2^n (n+1)\big)^{(n-1)c_w} \|P\|_w^{2(n-1)}
      |\xi-\alpha_t|_w^{t(t-1)}.
$$
Similarly, for all other places $v$ of $K$, we find with the same
definition of $c_v$
$$
|{\rm Disc}(P)|_v
  \le \big(2^n (n+1)\big)^{(n-1)c_v} \|P\|_v^{2(n-1)}.
$$
Applying the product formula we therefore obtain
$$
1 = \prod_v |{\rm Disc}(P)|_v
  \le  \big(2^n (n+1)\big)^{n-1} H(P)^{2(n-1)}
       |\xi-\alpha_t|_w^{t(t-1)}.
$$
The conclusion follows since we have $H(P)\ll
H(\alpha)^{n/\deg(\alpha)}$ for a root $\alpha$ of $P$ (see \S2).

\medskip
Our last result justifies the remark made in the introduction
concerning simultaneous approximation of several numbers by
conjugate algebraic numbers.

\proclaim Proposition 10.2. Assume that $K=\bQ$.  Let $n$ and $t$
be positive integers, and let $\kappa$ be a real number with
$$
\kappa > t^{-1}(t+1)^{1+(1/t)}.
$$
Then, there exist elements $\enum \xi t$ of $\bQ_w$ and a constant
$H_0 \ge 1$ (depending on $n$, $t$, $w$, $\kappa$) with the following
property. For any real number $H\ge H_0$ and any choice of numbers
$\enum \alpha t \in \Qbar_w$ that are algebraic over $\bQ$ of
degree $\le n$ and height $\le H$, we have
$$
\max_{1\le j\le t} |\xi_j-\alpha_j|_w \ge H^{-\kappa n^{1/t}}.
$$
\par

\medskip
Note that $t^{-1}(t+1)^{1+(1/t)}$ is a decreasing function of $t$
for $t>0$ tending to $1$ as $t$ tends to infinity.  For $t\ge 2$,
one can take $\kappa=3$.

\proof Put $b=(t+1)n$ and define a sequence of positive integers
$(a_\ell)_{\ell\ge 1}$ by the formula $a_\ell = [ b^{\ell/t} ]$.
Define also $\pi=1/2$ if $w$ is the place at infinity of $\bQ$ and
$\pi=p$ if $w$ corresponds to a prime number $p$.  We claim that
the elements of $\bQ_w$ given by
$$
\xi_j = \sum_{i=0}^\infty \pi^{a_{j+ti}},
  \quad
(j=1,\dots,t),
$$
have the required property.

To prove this, we choose a real number $\ep$ with $0<\ep<1$ such
that
$$
\kappa > {t+1+\ep \over t-\ep} (t+1)^{1/t},
$$
and consider the sequence of closed intervals $(I_\ell)_{\ell\ge
1}$ of $\bR$ given by
$$
I_\ell
  = \left[ \kappa^{-1} (t+1+\ep) a_\ell n^{1-(1/t)},\
           (t-\ep) a_\ell n \right].
$$
Two consecutive such intervals $I_\ell$ and $I_{\ell+1}$ overlap
for sufficiently large values of $\ell$ since
$$
\lim_{\ell\to\infty} {(t-\ep) a_\ell n \over
  \kappa^{-1} (t+1+\ep) a_{\ell+1} n^{1-(1/t)}}
  = {(t-\ep) \kappa \over (t+1+\ep) (t+1)^{1/t}}
  > 1.
$$
Therefore, the union of these intervals contains a half line
$[c,\infty)$ for some constant $c>0$.  Choose a real number $H$
with $H \ge |\pi|_w^{-c}$, and let $\enum \alpha t\in \Qbar_w$ be
algebraic over $\bQ$ of degree $\le n$ and height $\le H$.  By
definition of $c$, there exists an integer $\ell\ge 1$ such that
$-\log H/\log |\pi|_w\in I_\ell$. Writing $\ell$ in the form
$\ell=j+tm$ for some integers $j$ and $m$ with $1\le j\le t$ and
$m\ge 0$, we claim more precisely that, if $H$ is sufficiently
large (so that $\ell$ is large), we have
$$
|\xi_j-\alpha_j|_w \ge H^{-\kappa n^{1/t}}.
$$

For brevity, since the ultrametric case is simpler, we shall only
prove this refined claim in the case where $w=\infty$.  Then, we
have $\bQ_w=\bR$, $\Qbar_w=\bC$, $\pi=1/2$ and $\log H/\log 2\in
I_\ell$.  From now on, we drop the subscript $w$ on the absolute
value and consider the rational number
$$
r = \sum_{i=0}^m 2^{-a_{j+ti}}.
$$
Since $(a_{j+ti})_{i\ge 0}$ is a strictly increasing sequence of
positive integers, it satisfies
$$
\rmH(r)=2^{a_\ell}
  \and
2^{-a_{\ell+t}}
  \le \xi_j-r
   =  \sum_{i=1}^\infty 2^{-a_{\ell+ti}}
  \le 2^{-a_{\ell+t}+1}.
$$
If $\alpha_j=r$, we find
$$
|\xi_j-\alpha_j|
  = \xi_j-r
  \ge 2^{-a_{\ell+t}}
  \ge H^{-\kappa n^{1/t}}
$$
assuming that $H$ (and thus $\ell$) is sufficiently large so that
$$
a_{\ell+t} \le (t+1+\ep)na_\ell \le \kappa n^{1/t} \log H/\log 2.
$$
If $\alpha_j\neq r$, Liouville's inequality (see for example
Proposition 3.14 of [\refWa]) gives
$$
|\alpha_j-r|
  \ge \gamma \rmH(\alpha_j)^{-1} \rmH(r)^{-n}
  \ge \gamma H^{-1} 2^{-a_\ell n},
$$
with $\gamma=\gamma(n)=2^{1-n}(n+1)^{-1/2}$.  This implies
$$
|\alpha_j-r|
  \ge \gamma 2^{-(t-\ep)a_\ell n - a_\ell n}
  \ge 2^{-a_{\ell+t}+2}
  \ge 2 |\xi_j-r|,
$$
assuming that $H$ is sufficiently large so that
$$
(t+1-\ep)a_\ell n \le a_{\ell+t}-2+\log \gamma/\log 2.
$$
Since $\kappa >(t+1+\ep)/(t+\ep)$, we may also assume $(\gamma/2)
H^{-1} \ge H^{-((t+\ep)/(t+1+\ep))\kappa n^{1/t}}$ and so we get
$$
|\xi_j-\alpha_j|
  \ge |\alpha_j-r| - |\xi_j-r|
  \ge {1\over 2} |\alpha_j-r|
  \ge {\gamma\over 2} H^{-1} 2^{-a_\ell n}
  \ge H^{-\kappa n^{1/t}}.
$$

%
%

\section{References}

\narrower
\def\refmark#1{\medskip
\noindent \hskip -1.25 true cm \hbox to .9 true cm {\rm [#1]\hfill
} }

\refmark{\refBr} Brownawell, W.~D. -- Sequences of Diophantine
approximations. {\it J.~Number Theory} {\bf 6} (1974), 11--21.

\refmark{\refBT} Bugeaud, Y.; Teuli\'e, O. -- Approximation d'un
nombre r\'eel par des nombres alg\'ebriques de degr\'e donn\'e.
{\it Acta Arith.\ }{\bf 93}  (2000), 
77--86.

\refmark{\refBV} Bombieri, E.; Vaaler, J. -- On Siegel's lemma.
{\it Invent.\ Math.\ }{\bf 73} (1983), 11--32.

\refmark{\refARNAI} Davenport, H.; Schmidt, W.M. -- Approximation
to real numbers by algebraic integers. {\it Acta Arith.\ }{\bf 15}
(1969), 393-416.

\refmark{\refCartan} Cartan, H. -- {\it Cours de calcul
diff\'erentiel}. Hermann, Paris, 1977.

\refmark{\refCassels} Cassels, J.~W.~S. -- {\it An introduction to
the geometry of numbers}. Springer-Verlag, Berlin, 1971.

\refmark{\refCF} Cassels, J.~W.~S. -- Global fields. Chapter II
in: {\it Algebraic number theory}, J.~W.~S.~Cassels and
A.~Fr\"ohlich editors, Academic Press, 1967.

\refmark{\refFT} Fr\"ohlich, A.; Taylor, M. J. --- {\it Algebraic
number theory}. Cambridge Studies in Advanced Mathematics, {\bf
27}. Cambridge University Press, Cambridge, 1993.

\refmark{\refLa} Lang, S. -- {\it Fundamentals of Diophantine
geometry}. Springer-Verlag, New-York, 1983.

\refmark{\refLR} Laurent, M.; Roy, D. -- Criteria of algebraic
independence with multiplicities and interpolation determinants.
{\it Trans.\ Amer.\ Math.\ Soc.\ }{\bf 351} (1999), 1845-1870.

\refmark{\refMahler} Mahler, K. -- Ein \"Ubertragungsprinzip f\"ur
konvexe K\"orper. {\it \v{C}asopis P\v{e}st.\ Mat.\ }{\bf
68} (1939), 93--102.

\refmark{\refMahlerb} Mahler, K. -- Inequalities for ideal bases
in algebraic number fields. {\it J.\ Austral.\ Math.\ Soc.\ }{\bf
4} (1964), 425--448.

\refmark{\refMacfeat} Macfeat, R.~B. -- Geometry of numbers in
adele spaces.  {\it Dissertationes Math.\ (Rozprawy Mat.)} {\bf
88} (1971), 54 pp.

\refmark{\refRo} Roy, D. -- Approximation to real numbers by cubic
algebraic integers I. {\it Proc.\ London Math.\ Soc.} (to appear);
II. {\it Annals of Math.} (to appear).

\refmark{\refS} Schmidt, W.~M. -- On heights of algebraic
subspaces and diophantine approximations. {\it Annals of Math.\
}{\bf 85} (1967), 430--472.

\refmark{\refTe} Teuli\'e, O. -- Approximation d'un nombre
$p$-adique par des nombres alg\'ebriques. {\it Acta Arith.\ }{\bf
102} (2002), 137--155.

\refmark{\refWa} Waldschmidt, M. -- {\it Diophantine approximation on
linear algebraic groups: transcendence properties of the exponential
function in several variables}.  Springer-Verlag, New-York, 2000.

\refmark{\refWi} Wirsing, E. -- Approximation mit algebraischen
Zahlen beschr\"ankten Grades. {\it J.\ reine angew.\ Math.\ }{\bf
206} (1961), 67-77.

 \vskip 2truecm plus .5truecm minus .5truecm

\hbox{ \vtop{
 \hbox{Damien ROY}
 \hbox{D\'epartement de Math\'ematiques}
 \hbox{Universit\'e d'Ottawa}
 \hbox{585 King Edward}
 \hbox{Ottawa, Ontario K1N 6N5, Canada}
 \hbox{
      {
      droy@uottawa.ca}}
 \hbox{{
http{$:$}//aix1.uottawa.ca/${\scriptstyle \sim}$droy/}} }
\hfill\qquad \vtop{
 \hbox{Michel WALDSCHMIDT}
 \hbox{Universit\'e P.~et M.~Curie (Paris VI)}
 \hbox{Institut de Math\'ematiques CNRS UMR 7586}
 \hbox{Th\'eorie des Nombres, Case 247}
 \hbox{175, rue du Chevaleret, 75013 Paris, France}
 \hbox{
      {
      miw@math.jussieu.fr}}
 \hbox{{
http{$:$}//www.math.jussieu.fr/${\scriptstyle \sim}$miw/}} } }

\bye